\documentclass[times,12pt]{article}
\pdfoutput=1
\usepackage{amsmath,amssymb,amsfonts,latexsym,amsthm,enumerate,url}
\usepackage{graphicx}
\usepackage{mathrsfs}
\date{}

\newtheorem{thm}{Theorem}[section]

\newtheorem{lem}[thm]{Lemma}

\newtheorem{claim}{Claim}

\newcommand{\beq}[1]{\begin{equation}\label{#1}}
\newcommand{\enq}[0]{\end{equation}}

\newcommand{\remove}[1]{}

\newcommand{\comment}[1]{}

\title{On Shapley Values and Threshold Intervals}
\author{Gil Kalai\footnote{Einstein Institute of Mathematics, Hebrew University of Jerusalem and Efi Arazi School of Computer Science, Reichman University. Supported by ERC grant 834735 and by an ISF grant 2669/21.}~  and Noam Lifshitz\footnote{Einstein Institute of Mathematics, Hebrew University of Jerusalem. Supported by ERC grant 101163794 and by ISF grant  1980/22.}}

\begin {document}

\maketitle

\begin{abstract} 
Let $f\colon \{0,1\}^n\to \{0,1\}$ be a monotone Boolean functions, let $\psi_k(f)$ denote the Shapley value of the $k$th variable and $b_k(f)$ denote the Banzhaf value (influence) of the $k$th variable. We prove that if we have $\psi_k(f) \le t$ for all $k$, then the threshold interval of $f$ has length $\displaystyle O \left(\frac {1}{\log (1/t)}\right)$. We also prove that if $f$ is balanced and $b_k(f) \le t$ for every $k$, then
$\displaystyle \max_{k} \psi_k(f) \le O\left(\frac {\log \log (1/t)}{\log(1/t)}\right) $.
\end {abstract}


\section {Introduction}
\subsection {Boolean functions and threshold intervals}

Let $\Omega_n= \{0,1\}^n$ be the discrete $n$ dimensional cube. A function $f:\Omega_n \to \{0,1\}$ is called a 
Boolean function, and $f$ is \emph{monotone} if when $y_i \ge x_i$, for every $i$, then $f(x_1,x_2,\dots,x_n) \le f(y_1,y_2,\dots,y_n)$. For $p\in(0,1)$ the Boolean cube is equipped with the \emph{$p$-biased measure} $\mu_p,$ where the coordinates of $x\sim\mu_p$ are chosen independently with the value of each coordinate being one with probability $p$.  
We write 
$\mu_p(f)= \mathbb{E}_{x\sim \mu_p}[f(x)].$
For $\epsilon >0$ we write $p_\epsilon (f)$ as the value $q$ such that $\mu_q(f)=\epsilon$. Given a small real number $\epsilon$ we refer to the interval $[p_\epsilon (f),p_{1-\epsilon}(f)]$ as the {\it $\epsilon$-threshold interval} of $f$. For instance, when $f$ is the majority function, the length of its $\epsilon$-threshold interval is $C(\epsilon) \frac {1}{\sqrt n}$.  

\subsection {Cooperative games and the Shapley value}

Cooperative $n$-player games (with transferable utilities) are defined via a general payoff function $f$ that assigns for a subset $S \subset [n]$  a real number $f(S)$, its payoff $f(S)$. We will assume that $f(\emptyset)=0$.

Boolean functions can be regarded as special instances of a cooperative voting game.  
Every vector $x=(x_1,x_2,\dots,x_n) \in \Omega_n$ is associated to a coalition $S_x=\{i: x_i=1\}$, and we regard $f$ as a payoff function for coalitions $S_x$. A coalition $S_x$ is a \emph{winning coalition} if $f(x)=1$.

The {\it Shapley value} \cite {Sha53} is a central solution concept for cooperative games based on a system of axioms. For simple games it is also referred to as the {\it Shapley--Shubik power index} \cite {ShaShu54}. 
For $k$ with $1\le k \le n$ and an ordering $(\pi(1),\pi(2),\dots, \pi (n))$ of $[n]$, let $s$ be such that $k=\pi (s)$.
Consider the subset of $[n]$, $S(\pi,k)=\{\pi(1),\pi(2),\dots \pi(s-1)$ and 
Shapley's original definition was given by the following formula
\begin {equation}
\psi_k(f) = \frac {1}{n!}\sum_{\pi} f(S(\pi,k) \cup \{k\}) - f(S(\pi,k).    
\end {equation}
In words, we consider the average contribution of player $k$ to the preceeding players taken over all $n!$ orderings of the players. 

Shapley proved that the vector ($\psi_1(f), \psi_2(f),\dots \psi_n(f))$ satisfies and is uniquely determined by the following axioms:
\begin {enumerate}
\item  {\it Efficiency.} $$\sum \psi_k(f) = f([n]).$$
\item  {\it Symmetry.} Let $i\ne j\in [n]$. Suppose that for every $S\subseteq [n]$ with $i,j\notin S$ we have
$$f(S \cup \{i\})=f(S \cup \{j\}).$$
Then $$\psi_i(f)=\psi_j(f).$$
\item  {\it Null-player.} Suppose that  $k \in [n]$ and for every $S \subset [n]$ 
$$ f(S \cup \{k\})=f(S).$$
Then $$\psi_k(f)=0.$$
\item  {\it Additivity.} For two payoff-functions $f$ and $g$
$$\psi_k(f+g)=\psi_k(f)+\psi_k(g).$$
\end {enumerate}

Moreover, the Shapley value satisfies also the additional property that for every real number $a$ and 
a vector $b= (b_1,b_2,\dots,b_n) \in \mathbb R^n$, if $g(S)$ is defined by $g(S)=af(S)+\sum_{i \in S}b_i$, then 
$\psi_i(g)=a\psi_i(f)+b_i.$
For a survey on the Shapley value, alternative axiomatic approaches, and some applications to real-life voting rules see Winter \cite {Win02}. Recently, the Shapley value has served as the foundation for multiple methods that attribute a machine-learning model's prediction for a given input to its individual features (see, e.g.  \cite {LunLee17}).

\subsection{Influence and the Banzhaf value}

Let
$$\sigma_k(x_1,\dots,x_{k-1},x_k,x_{k+1},\dots,x_n)
= (x_1,\dots,x_{k-1},1-x_k,x_{k+1},\dots,x_n).$$

The influence $I^p_k(f)$ of the $k$th variable on a Boolean function $f$ is defined by:
$$I^p_k(f) = \mu_p(x \in \Omega_n, f(x) \ne f(\sigma_k(x))).$$
The total influence is defined by $$I^p(f)=\sum_{k=1}^nI^p_k(f).$$

The notion of total influence is key in the study of threshold intervals is key due to 
Russo's lemma \cite {Rus82}, which asserts that 
$$\frac {d\mu_p(f)}{dp}=I^p(f).$$



The {\it Banzhaf value} \cite {Ban64} of the $k$-the player $b_i(f)$ is its influence $I^p_k(f)$, for $p=1/2$. 
Owen \cite {Owe88} found the following alternative definition for Shapley value  

$$\psi_k(f) = \int_{0}^{1} I_k^p(f)dp.$$

Note that Owen's theorem and Russo's Lemma gives together a quick verification of the efficiency property.


\subsection{The KKL theorem, sharp threshold, and our previos work}

The KKL theorem \cite {KKL88} asserts that if $f$ is a Boolean function and $\mu(f)=t$ then 
$$\max_{1 \le k \le n} I_k(f) \ge Ct(1-t)\frac {\log n}{n}.$$ 

KKL's theorem implies that if all influences are the same then the sum of influences is at least $C t (1-t)\log n $, and some further results from \cite {KKL88} shows that for this conclusion to hold it suffices that $\max_{1 \le k \le n} I_k(f) < n^{-c}$ for some $c>0$. These results were extended to $\mu_p$ and were used to show a ``sharp thershold phenomena" in several subsequent papers \cite {BKKKL92,Tal94,FriKal96}. 

Let $(f_n)$ be a sequence of monotone balanced Boolean functions. We assume for simplicity that $f_n$ is defined on $\Omega_n$. The sequence $(f_n)$ has a sharp threshold behavior if for every fixed $\epsilon$, $0< \epsilon<1/2,$ 
$$\lim_{n \to \infty}(p_{1-\epsilon}(f_n)- p_\epsilon (f_n)) =0.$$
In \cite {Kal04}, Kalai proved that for a sequence $(f_n)$ of monotone balanced Boolean functions a necessary and sufficient condition for the sharp threshold behavior, is that
$$ \lim_{n \to \infty}\max _{1 \le k \le n}\psi_k(f_n)=0.$$ In this paper we give a sharp form of this result.  

\subsection {Our main results}
Our main result is the following: 

\begin {thm}\label{thm:main theorem}
There exists an absolute constant $C>0$, such that the following holds. Let $f$ be a monotone Boolean function, $\epsilon,t>0$. Suppose that $p_{1-\epsilon}(f)-p_\epsilon (f) >t$. 
Then there exists a set $T\subset [n]$ of size $\le \epsilon^{-Ct}$ such that 
for every $p$ with $ p_\epsilon \le p \le p_{1-\epsilon}$ there exist $i \in T$ such that 
$$I_i^p(f) \ge \epsilon^{-Ct}.$$ 
\end {thm}

This theorem implies the following two results

\begin {thm}
\label {t:main}
Let $f$ be a monotone Boolean function. If for all $k\in [n]$ we have $\psi_k(f) \le t$, then 
$$p_{1-\epsilon}(f)-p_\epsilon (f) \le C \frac{\log (1/\epsilon)}{\log (1/t)}.$$
\end {thm}

A weaker form of Theorem \ref {t:main} was proved in \cite {Kal04} and the argument there yields only 
that $\displaystyle p_{1-\epsilon}(f)-p_\epsilon (f) \le C \frac{\log (1/\epsilon)}{\log \log (1/t)}.$

The dependency on $t$ in Theorem \ref {t:main} is sharp due to the tribes example described below, which is an example with all Shapley values equal to $1/n$ and with an $\epsilon$-threshold interval of length $\ge \frac{1}{C\log n}$ for $\epsilon\le 1/3$ and an absolute constant $C>0$.
\begin {thm} 
\label {t:shaba}
Let $f$ be a monotone Boolean function and let  $t,\epsilon >0$. Suppose 
that $\epsilon \le \mu_{1/2}(f) \le 1-\epsilon$. Then if 
$$ b_k(f) \le t$$ for all $k$, then  we have 
$$\max_k \psi_k(f) \le \max\left(\frac{C\log\log (1/t)}{\log(1/t)}, \frac{C\log(1/\epsilon)}{\log(1/t)} \right).$$
\end {thm}

Theorem \ref {t:shaba} implies that if a sequence of strong simple games $f_n$ 
on $n$ variables with $\mu_{1/2}(f_n)\in (0.1,0.9)$ satisfies that 
$\lim_ {n \to \infty} \max_{1\le k\le n} b_k(f) = 0$, then we also have $\lim _{n \to \infty} \max_{1\le k\le n} \psi_k(f) = 0.$ The reverse implication does not hold as seen from the first  example in the next subsection. The second example below shows that Theorem \ref{t:shaba} is nearly optimal in the sense that there exists a Boolean function on $n$ variables with $\max_k b_k(f) = t = O(\log n/n)$ while $\max_k \psi_k(f)\ge c/\log(1/t) $ for an absolute constant $c>0$.

\subsection{Two examples} 
\subsubsection *{Supermajority or else a judge} 
Consider the voting rule between two candidates where there are $n$ voters and 
the candidate who received the majority of votes wins the election unless the gap between the votes for the two candidates is smaller than $n^{2/3}$ in which case $f(x_1,x_2,\dots,x_n)=x_1$.
(We can think about the first voter as the ``judge".) In this example the first player's Banzhaf index is close to one. In contrast, all Shapley values are  $O(n^{-1/3})$. 
\subsubsection *{Supermajority and a judge or else tribes}
Let $g(x_1,x_2,\dots,x_n)$ be the Ben-Or--Linial tribe example defined as follows: divide $[n]$ into tribes $J_1,J_2,\dots,J_m$ of size $\log n -\log \log n + \log \log e$ and then define $g(x_1,x_2,\dots,x_n)=1$ if for some tribe $J_k$, $x_i=1$ for every $i \in J_k$. The parameters are chosen for this function to satisfy $\mu_{1/2}(f)= 1/2(1+o(1))$. Next, 
let $  \delta = \frac{1}{1000\log n}$, and  consider the Boolean function $f(x_1,x_2,\dots,x_n)$, where $f(x)=1$ if either 
a) $x_1=1$ and $|\{i: x_i=1\} > \frac {1}{2}+ \delta$,
or $g(x_1,x_2,\dots ,x_n)=1$. For this example, the Banzhaf values of all players are roughly equal to $O(\frac{\log n}{n})$, while the maximum Shapley value (attained for $i=1$) is $\ge \frac {1}{C\log n}$ for an absolute constant $C>0$.   

\subsection {Voting paradoxes, Condorcet and McGarvey}

Consider a monotone voting rule based on a Boolean function $f(x_1,x_2,\dots,x_n)$ satisfying $$f(1-x_1,1-x_2,\dots, 1-x_n)=1-f(x_1,x_2,\dots,x_n).$$ When there are $m$ candidates and $n$ voters,  each with an ordered preference relation over the candidates, $f$ produces an asymmetric social preference relation. McGarvey \cite{McG53} showed that the majority rule can yield any asymmetric social preference relation when the number of voters is large relative to the number of candidates. This result generalizes the well-known {\it Condorcet Paradox} for three candidates. 

Theorem \ref {t:main} together with an argument from Kalai \cite {Kal04} and a result of Alon \cite {Alo02}, implies that for some absolute constant $C$, when there are $m$ candidates, a voting rule based on $f$, for which the Shapley values of every voter is at most $e^{-C\sqrt {m} \log m}$, leads to every asymmetric preference relation. This provides a quantitative version of a result from \cite {Kal04}; however, its sharpness remains an open question. 

The sharp threshold behavior of the majority function follows from the weak law of large numbers. {\it Condorcet's Jury Theorem} provided an early political science interpretation; see, e.g., \cite{AusBan96, Mos22}. Theorem \ref{t:main} establishes a sharp connection, for general voting rules, between the assertion of Condorcet's Jury Theorem and the vanishing of the maximum Shapley value.



\section{Preliminaries}
\subsection {$p$-biased characters}
Let $p\in  (0, 1)$ be a real number. 
The space $L_2(\mu_p)$ is the space of real functions on $\Omega_n$ with the inner product 
$$\langle f,g\rangle =\mathbb E _{x \sim \mu_p} f(x)g(x).$$
For $S \subset [n]$ we let $$\chi_S^p(x)=\prod_{i \in S} \frac {x_i-p}{\sqrt {p(1-p)}}.$$
These {\it characters} form an orthonormal basis for $L_2(\mu_p)$. 

\subsection{One sided noise operator}

The one sided noise operator $\mathrm{T}^{p\to q}$ is defined on $\Omega_n$ as follows. Given $x \in \Omega_n$
$N(x)$ is a distribution on $\Omega_n$ where every coordinate with value 0 stays the same and every coordinate with value 1 becomes a 0 with probability $\frac {q-p}{q}$, independently over all $i$. The one sided noise operator is then defined via:
$$\mathrm{T}^{p\to q}(f)(x)= \mathbb E _{y\sim N(x)}f(y).$$

Note that if $x \sim \mu_q$ and $y \sim N(x)$, then the marginal distribution for $y$ is $\mu_p$. 

If $f$ is monotone then $$\langle \mathrm{T}^{p\to q},f\rangle_{\mu(q)}=\mathbb E_{x \sim \mu_q, y \sim N(x)}f(x)f(y) =\mu_p(f).$$ Our argument will be based on finding conditions implying that 
$f(x)$ and $f(y)$ are nearly independent, and therefore the above expectation should be actually close to $\mu_p(f)\mu_q(f)$. 
This will allow us to deduce an upper bound on the size of the threshold result because it implies that for every $p,q$ that are far apart $\mu_p(f)$ is close to $\mu_p(f)\mu_q(f)$, in which case either $\mu_p(f)$ is close to $0$ and $p$ is not in the threshold interval or $\mu_q(f)$ is close to $1$ and $q$ is not.

\subsection*{Fourier formula for the one sided noise operator}
In turns out that we have the following convenient Fourier formula for the one sided noise operator.
\begin{lem}[\cite {Lif20}]
\label{lem:Eigenvalues of one sided noise operator} Let $p<q\in\left(0,1\right)$,
$\rho=\sqrt{\frac{p\left(1-q\right)}{q\left(1-p\right)}}.$ Let $f\in L^{2}\left(\left\{ 0,1\right\} ^{n},\mu_{p}\right)$
and write 
\[
f=\sum\hat{f}\left(S\right)\chi_{S}^{p}.
\]
 Then 
\[
\mathrm{T}^{p\to q}f=\sum\rho^{\left|S\right|}\hat{f}\left(S\right)\chi_{S}^{q}
\]
\end{lem}
\subsection {A KKL-type argument} 

We write $D_{i}\left[f\right]=\sum_{S\ni i}\hat{f}\left(S\right)\chi_{S\setminus\left\{ i\right\} },$
and $\mathrm{Inf}_i[f] = \|D_i[f]\|_2^2.$ Given a function $f$ on the Boolean cube $i\in [n]$ and $b\in \{0,1\}$ let us write $f_{i\to b}$ for the function on $\{0,1\}^{[n]\setminus\{i\}}$ with $f_{i\to b}(x) = f(x,b),$ where $b$ is plugged in the $i$th coordinate. Then we have $D_i[f](x) = \sqrt{p(1-p)}(f_{i\to 1} - f_{i\to 0})$, an equality of linear operators that can be easily verified on the basis of characters. 
We have $I_i^p[f] = \|f_{i\to 1} - f_{i\to 0}\|_2^2$ and so $\mathrm{Inf}_i^p[f] = p(1-p) I_i^p[f].$ The \emph{degree} of a function is the maximal $d$ for which there exists $S$ with $\hat{f}(S)\ne 0.$ When $p$ is clear from context we write $f^{\le d}$ for the truncation given by  \[f^{\le d}  = \sum_{|S| \le d} \hat{f}(S)\chi_S^{p}.\] We shall also write $f^{=d}$ for $f^{\le d} - f^{\le d-1}.$

We shall make use of the following consequence of the hypercontractivity theorem over the $p$-biased cube.
\begin{thm}[{\cite[Chapter 10]{Odo14}}]\label{thm:biased bonami}
    Let $f$ be a function of degree $d$, and let $\lambda = \min(p,1-p).$ Then 
    \[
    \|f\|_{L_4(\mu_p)} \le 3^{d/2}\lambda^{-d/4}\|f\|_{L_2(\mu_p)}.
    \]
\end{thm}

As a consequence of Theorem \ref{thm:biased bonami} we have the following bound on the low degree part of the derivatives $D_i[f]$.

\begin{lem}
\label{lem:KKL}Let $d>0$, $\lambda = \min(p,1-p)$. Then we have 
\[
\|D_{i}\left[f\right]^{\le d}\|_{L_2(\mu_{p})}^{2}\le 3^{d} \lambda^{-d/2}\sqrt{p(1-p)} I^{p}_i[f]^{1.5}.
\]
\end{lem}
\begin{proof}
By H\"{o}lder's inequality and Theorem \ref{thm:biased bonami} we have
\begin{align*}
\|D_{i}\left[f\right]^{\le d}\|_{2}^{2} & =\left\langle D_{i}\left[f\right]^{ \le d},D_{i}\left[f\right]\right\rangle \\
 & \le\|D_{i}\left[f\right]^{\le d}\|_{4} \|D_{i}\left[f\right]\|_{\frac{4}{3}}\\
 & \le \lambda^{-d/4}3^{d/2}\|D_{i}[f]^{\le d}\|_{2}\|D_{i}\left[f\right]\|_{\frac{4}{3}}.
\end{align*}
 Rearranging, we obtain 
\[
\|D_{i}\left[f\right]^{\le d}\|_{2} \le \lambda^{-d/4} 3^{d/2} \|D_{i}\left[f\right]\|_{\frac{4}{3}}
\]
 and the lemma follows from the fact that \[\|D_{i}\left[f\right]\|_{4/3} = \sqrt{p(1-p)} \|f_{i\to 1} - f_{i\to 0}\|_{4/3} = \sqrt{p(1-p)} I_i[f]^{3/4}.\]
\end{proof}

\section{Proof of results}

\begin{lem}
\label{lem:Main Lemma} Let $p<q, \rho = \sqrt{\frac{p(1-q)}{q(1-p)}}, \lambda = \min(p,1-p), r= \max(0, -\log_{\rho}(6\rho \lambda^{-1/4})).$
Let $f,g\colon\left\{ 0,1\right\}^{n} \to \left\{ 0,1\right\}$.
Let $\epsilon >0$ and suppose that \[ \max_i (I_i^{p}[f] \cdot \mathrm{Inf}_i^{q}[g]) \le I^{p}[f]^{-5} (\epsilon/2)^{8 + 8r}.\]

Then
\[
\left\langle \mathrm{T}^{p\to q}f,g\right\rangle \le\mu_{p}\left(f\right)\mu_{q}\left(g\right)+\epsilon.
\]
\end{lem}

\begin{proof}

Let $d=\left\lfloor \log_{\rho}\left(\epsilon/2\right)\right\rfloor ,$
and let 
$\tau_{1},\text{\ensuremath{\tau_{2}} }$ to be chosen later
with 
\[
\min\left\{ \tau_{1},\tau_{2}\right\} >\frac{\epsilon^{C\left(As+1\right)}}{I^{p}\left[f\right]^{20}},
\]
 where $C$ is a sufficiently large absolute constant. Write 
\[
T_{1}=\left\{ i:\,\mathrm{Inf}_{i}^{p}\left[f\right]\ge\tau_{1}\right\} ,
\]
 and 
\[
T_{2}=\left\{ i:\mathrm{Inf}_{i}^{q}\left[g\right]\ge\tau_{2}\right\} .
\]
 Then by hypothesis $T_{1}\cap T_{2}=\varnothing.$ Write $f=\sum\hat{f}\left(S\right)\chi_{S}^{p},$
and $g=\sum\hat{g}\left(S\right)\chi_{S}^{q}.$

By Lemma \ref{lem:Eigenvalues of one sided noise operator}, and Cauchy--Schwarz
we have 
\begin{align*}
\left\langle \mathrm{T}^{p\to q}f,g\right\rangle  & =\sum_{S\subseteq\left[n\right]}\rho^{\left|S\right|}\hat{f}\left(S\right)\hat{g}\left(S\right)\\
 & \le\frac{\epsilon}{2}\sqrt{\sum_{\left|S\right| > d}\hat{f}\left(S\right)^{2}}\sqrt{\sum_{\left|S\right| > d}\hat{g}\left(S\right)^{2}}+\\
 & +\sum_{k=0}^{d} \rho^{k}\sum_{\left|S\right| = k}\hat{f}\left(S\right)\hat{g}\left(S\right)\le \epsilon/2 +\sum_{k=0}^{d} \rho^{k}\sum_{\left|S\right| = k}\hat{f}\left(S\right)\hat{g}\left(S\right).
\end{align*}
 
By Cauchy--Schwarz and Lemma \ref{lem:KKL} we have  
\begin{align*}
\left| \sum_{|S| = k }\hat{f}\left(S\right)\hat{g}\left(S\right) \right|& = \frac{1}{k} \left| \sum_{i=1}^{n}\sum_{|S| = k,i\in S}\hat{f}\left(S\right)\hat{g}\left(S\right) \right|\\
 & \le \frac{1}{k}\sum_{i=1}^{n}\sqrt{\sum_{S\ni i,\left|S\right| = k}\hat{f}\left(S\right)^{2}}\sqrt{\sum_{S\ni i,\left|S\right| = k }\hat{g}\left(S\right)^{2}}\\
 & = \frac{1}{k}\sum_{i=1}^{n}\|D_{i}\left[f\right]^{=k-1}\|_{L_2(\mu_p)}\|D_{i}\left[g\right]^{=k-1}\|_{L_2(\mu_q)}\\
 & \le \frac{\lambda^{-(k-1)/4}3^{(k-1)/2}}{k}p^{1/4}(1-p)^{1/4} I_i^{p}[f]^{3/4} \|D_{i}\left[g\right]^{=k-1}\|_{L_2(\mu_q)}.
\end{align*}

Let $\tau = \max_i(I^p_i[f] \mathrm{Inf}^p_i[g]).$ 
Then we may upper bound 
\begin{align*}
    I_i^{p}[f]^{3/4} \|D_{i}\left[g\right]^{=k-1}\|_{L_2(\mu_q)} & = I_i^p[f]^{5/8} \|D_{i}\left[g\right]^{=k-1}\|_{L_2(\mu_q)}^{3/4} I_i^p[f]^{1/8} \|D_{i}\left[g\right]^{=k-1}\|_{L_2(\mu_q)}^{1/4} 
    \\ & \le I_i^p[f]^{5/8} \|D_{i}\left[g\right]^{=k-1}\|_{L_2(\mu_q)}^{3/4} \tau^{1/8}.
\end{align*}
Therefore, summing over all $i$ and applying H\"{o}lder's inequality we have 
\begin{align*}
 \sum_{i=1}^n  I_i^{p}[f]^{3/4} \|D_{i}\left[g\right]^{=k-1}\|_{L_2(\mu_q)} & \le \tau^{1/8} I^p[f]^{5/8}\left(\sum_{i=1}^{n} \|D_i[g]^{=k-1}\|_2^2\right)^{3/8} \\ &= \tau^{1/8} I^p[f]^{5/8} k^{3/8}\|f^{=k}\|_2^{3/4} \\ & \le \tau^{1/8} I^p[f]^{5/8} k^{3/8}. 
\end{align*}

Putting everything together yields 

\begin{align*}
\left\langle \mathrm{T}^{p\to q}f,g\right\rangle -\mu_p(f)\mu_p(g) & \le \epsilon/2 + \tau^{1/8} \sum_{k=1}^{d} 2^{-k} (6\rho \lambda^{-1/4})^{k} I^{p}[f]^{5/8} \\& \le \epsilon/2 + (\epsilon/2)^{-r} \tau^{1/8} I^{p}[f]^{5/8} \le \epsilon.   
\end{align*}

\end{proof}
Applying the above lemmas on the functions given by $f^{\dagger}\left(x\right)=f\left(\overline{x}\right),g^{\dagger}\left(x\right)=g\left(\overline{x}\right)$
immediately yields the following statement with the roles of $p,q$ reversed.
\begin{lem}
\label{lem:dual lemma} Let $p<q, \rho = \sqrt{\frac{p(1-q)}{q(1-p)}}, \lambda = \min(q,1-q), r= \max(0, -\log_{\rho}(6\rho\lambda^{-1/4})).$
Let $f,g\colon\left\{ 0,1\right\}^{n} \to \left\{ 0,1\right\}$.
Let $\epsilon >0$ and suppose that 
\[
\max_i (I_i^{p}[f] \cdot \mathrm{Inf}_i^{q}[g]) \le I^{q}[g]^{-5} (\epsilon/2)^{8 + 8r}.
\]
Then
\[
\left\langle \mathrm{T}^{p\to q}f,g\right\rangle \le\mu_{p}\left(f\right)\mu_{q}\left(g\right)+\epsilon.
\]
\end{lem}

Finally, we have the following simplified variant of Lemmas \ref{lem:Main Lemma} and \ref{lem:dual lemma}.

\begin{lem}\label{lem:simplified main lemma}
    There exists an absolute constant $C>0$, such that the following holds. Let $p<q,\epsilon \in (0,1/2)$ and $f,g\colon\{0,1\}^n\to \{0,1\}$ be with 
    \[
    \max_i(I_i^p[f]I_i^p[g]) \le \min(I^p[f]I^q[g])^{-5} \epsilon^{\frac{C}{q-p}}.
    \]
    Then $\left\langle T^{p\to q}f,g\right\rangle \le \mu_p(f)\mu_q(g) +\epsilon.$
\end{lem}
\begin{proof}
    The statement follows immediately from either Lemma \ref{lem:Main Lemma} or Lemma \ref{lem:dual lemma} noting that we have the crude upper bound  $\log_{\rho} (6\rho\lambda^{-1/4}) \le  O(\frac{1}{q-p}).$
\end{proof}

If $f\colon\left\{ 0,1\right\} ^{n}\to\left\{ 0,1\right\} $ is a
non-constant monotone Boolean function, then we write $p_{\alpha}=p_{\alpha}\left(f\right)$
for the unique value of $p$ for which $\mu_{p_{\alpha}}\left(f\right)=\alpha.$
\begin{thm}\label{thm:main theorem again}
There exists $C>0$, such that the following holds. Let $f$ be a
monotone Boolean function, $q=p_{1-\epsilon},p=p_{\epsilon}$. Let
\[
\epsilon\in(0,1/2),s=\frac{1}{q-p},
\]
 
Then there exists a set $T$ of size $\le\epsilon^{-Cs},$
such that for all $r\in\left[p,q\right]$ there exists an $i\in T$
with $\mathrm{Inf}_{i}^{\left(r\right)}\left[f\right]>\epsilon^{Cs}.$
\end{thm}

\begin{proof}
By Russo's lemma we may choose $p'\in\left[p+\frac{p-q}{6},p+\frac{p-q}{3}\right]$
with $I^{p'}\left[f\right]\le 6s.$ Let $C_{1}$ be
a sufficiently large constant, and let 
\(
\tau_{1}= \epsilon^{C_{1} s}.
\)
Set $T_{1}=\left\{ i:\mathrm{Inf}_{i}^{\left(p'\right)}\left[f\right]\ge\tau_{1}\right\} ,$
then $\left|T_{1}\right|\le\epsilon^{-C_1s}6s \le \epsilon^{-Cs},$ provided that $C$ is sufficiently large. 
Then
by Lemma \ref{lem:simplified main lemma} (with $p',\tilde{q},\frac{\epsilon^{2}}{2}$)
for each $\tilde{q}\in\left[\frac{p+q}{2},q\right]$ there exists
some $i\in T_{1}$ with $\mathrm{Inf}_{i}^{\left(\tilde{q}\right)}\left[f\right]\ge\tau_{1}.$
Indeed, otherwise we would have 
\[
\mu_{p'}\left(f\right)=\left\langle \mathrm{T}_{p'\to\tilde{q}}f,f\right\rangle \le\mu_{p'}\left(f\right)\mu_{\tilde{q}}\left(f\right)+\frac{\epsilon^{2}}{2},
\]
 which contradicts $\mu_{p'}\left(f\right)\left(1-\mu_{\tilde{q}}\left(f\right)\right)\ge\epsilon^{2}.$

Similarly letting $q'\in\left[q-\frac{p+q}{3},q-\frac{p+q}{6}\right]$
with $I^{q'}\left[f\right]\le6s$ 
\(
\tau_{2}= \epsilon^{C_{1}s}
\)
 and $T_{2}=\left\{ i:\mathrm{Inf}_{i}^{\left(p'\right)}\left[f\right]\ge\tau_{2}\right\} $
we obtain that $\tau_{2}\ge\epsilon^{Cs},\left|T_{2}\right|\le\epsilon^{-Cs}$
and for each $\tilde{p}\in\left[p,\frac{p+q}{2}\right]$ we have $\mathrm{Inf}_{i}^{\left(\tilde{p}\right)}\left[f\right]\ge\tau_{2}$
for some $i\in T_{2}.$ Setting $T=T_{1}\cup T_{2}$ completes the
proof.
\end{proof}
The result about the Shapely value follows easily from the theorem.

Theorem \ref{thm:main theorem} is the same as Theorem \ref{thm:main theorem again} stated twice without a good reason.  

\begin{proof}[Proof of Theorem \ref{t:main}]
   Let $s = \frac{1}{p_{1-\epsilon} -p_{\epsilon}}.$ By Theorem $\ref{thm:main theorem}$ there exists an absolute constant $C>0$ and a set $T$ of size $\le \epsilon^{Cs}$, such that for all $p\in [p_{\epsilon}, p_{1-\epsilon}]$ we have $\sum_{i\in T}\mathrm{I}_i^{(p)}[f] \ge \epsilon^{Cs}$. Integrating over the $\epsilon$-threshold interval, which has size $1/s$, this implies that $\sum_{i\in T}\psi_i(f) \ge 1/s\epsilon^{Cs}\ge \epsilon^{2Cs}$.This shows that for some $i\in T$ we have $\psi_i(f) \ge \epsilon^{2Cs}/|T| \ge \epsilon^{3Cs}.$ Therefore if $\psi_i(f)\le t$ for all $i$, then $s \ge \frac{\log(1/t)}{3C \log(1/\epsilon))}$. 
\end{proof}

\begin{proof}[Proof of Theorem \ref{t:shaba}]
If $\epsilon \ge 1/\log(1/t),$ then the hypothesis still holds with $\epsilon$ replaced by $1/\log(1/t)$ and the conclusion stays the same. Therefore we may assume without loss of generality that $\epsilon \le \frac{1}{\log(1/t)}.$
We may upper bound $\mathrm{I}_i^p[f]$ by $1$ inside the threshold interval to obtain that $\int_{p_{\epsilon}}^{p_{1-\epsilon}}I_i^{p}[f]dp \le p_{1-\epsilon} - p_{\epsilon}.$ By Russo's lemma we may also upper bound 
\[
\int_{0}^{p_{\epsilon}}\mathrm{I}_i^p[f]dp \le \int_{0}^{p_{\epsilon}}\sum_{i=1}^n\mathrm{I}_i^p[f]dp =\mu_{p_{\epsilon}}(f) = \epsilon
\]
and similarly,
\[
\int_{p_{1-\epsilon}}^{1}\mathrm{I}_i^p[f]dp \le \epsilon. 
\]
This shows that 
$\psi_i(f) \le 2\epsilon + p_{1-\epsilon} -p_{\epsilon}.$ 
Now by Theorem \ref{thm:main theorem} with $p=1/2$ there exists $i$ with $b_i= I^{1/2}_i[f] \ge \epsilon^{-Cs},$ where $s =\frac{1}{p_{1-\epsilon} - p_{\epsilon}}.$ As earlier rearranging yields that $p_{1-\epsilon} -p_{\epsilon} = O(\frac{\log(1/\epsilon)}{\log(1/t)}).$
Now $2\epsilon + O(\frac{\log(1/\epsilon)}{\log(1/t)}) \le   O(\frac{\log(1/\epsilon)}{\log(1/t)}),$ since $\epsilon \le 1/\log(1/t).$ This completes the proof of the statement.
\end{proof}


\begin{thebibliography}{99}

\bibitem {Alo02} N. Alon., Voting Paradoxes and Digraphs Realizations,{\it Advances in Applied Mathematics}, 29 (2002), 126--135.

\bibitem {Arr50} K. Arrow, A difficulty
in the theory of social welfare, {\it Journal of Political Economy} 58 (1950) , 328--346.

\bibitem {AusBan96} D. Austen-Smith and JS Banks, Information aggregation, rationality, and the
 Condorcet Jury Theorem, {\it American Political Science Review} 90(1996)34–45

\bibitem {Ban64}
J. F. Banzhaf,  Weighted voting doesn’t work: A mathematical analysis. Rutgers Law
 Review, 19 (1964), 317.
 
\bibitem {BenLin85}
M. Ben-Or and N. Linial, Collective coin flipping, robust voting games,
and minima of Banzhaf value, {\it Proc. 26th IEEE Symp. on the
Foundation of Computer Science} (1985), pp. 408--416.
\bibitem {BenLin90} M. Ben-Or and N. Linial,
Collective
coin flipping, in {\it Randomness and Computation} (S. Micali, ed.),
New York, Academic Press (1990), pp. 91--115. 


\bibitem {BolTho87} B. Bollob\'as and A. Thomason, Threshold functions,
 {\it Combinatorica} 7 (1987), 35--38.
 
\bibitem{BKKKL92} J. Bourgain, J. Kahn, G. Kalai, Y. Katznelson
and N. Linial,
The influence of variables in product spaces, {\it Israel J. Math.} { 77}
(1992), 55--64.

 \bibitem{Fri98} E. Friedgut, Boolean functions with low average
sensitivity depend on few coordinates, {\it Combinatorica} 18 (1998), 27--35.


\bibitem {FriKal96} E. Friedgut and G. Kalai,
Every monotone graph property has a sharp threshold,
{\it Proc.\ Amer.\ Math.\ Soc.} {\bf 124} (1996), 2993--3002.


\bibitem {KKL88}
J. Kahn, G. Kalai, and N. Linial, The influence of variables
on Boolean functions, in {\it Proc. 29-th Annual Symposium on Foundations of
Computer Science} (1988), 68--80.

\bibitem {Kal04} G. Kalai, Social indeterminacy, {\it Econometrica,} 72 (2004), 1565--1581.
\bibitem{Led01}
M. Ledoux, {\it The Concentration of Measure Phenomenon,}
Mathematical Surveys and Monographs, 89,
American Mathematical Society, Providence, RI (2001).

\bibitem {LunLee17} S. M. Lundberg and S.-I. Lee, A unified approach to interpreting model predictions, in {\it Proc. 31st Int. Conf. Neural Information Processing Systems} pp.4768–4777 (2017).

\bibitem {Lif20}  N. Lifshitz, Hypergraph removal lemmas via robust
 sharp threshold theorems, Discrete Analysis 2020:10, 46 pp, arXiv:1804.00328.
\bibitem {McG53}

D. C. McGarvey, A theorem on the construction of voting paradoxes, Econometrica
21 (1953), 608--610.

\bibitem {Mos22} E. Mossel, Probabilistic view of voting, paradoxes, and manipulation, {\it Bull. Amer. Math. Soc.} 59 (2022), 297--330.

\bibitem {Odo14} R. O'Donnell, {\it Analysis of Boolean Functions,}
Cambridge University Press, 2014.


\bibitem {Owe88} G. Owen, 
Multilinear extensions of games. {\it The Shapley value},
139--151, Cambridge Univ. Press, Cambridge, 1988.
 
\bibitem {Sha53}
L. S. Shapley, A value for n-person games. In Contributions to the Theory of 
Games, vol II, ed. H.W. Kuhn and A.W. Tucker. Princeton: Princeton University Press, 1953.


 \bibitem {ShaShu54}
 L. S. Shapley and M. Shubik,  A method for evaluating the distribution of power in a
 committee system. American Political Science Review, 48 (1954), 787–792.


\bibitem {Rus82}
L. Russo, An approximate zero-one law,
{\it  Zeitschrift f\"ur Wahrscheinlichkeitstheorie und Verwandte Gebiete}
61 (1982), 129--139.

\bibitem {Tal94}
M. Talagrand, On Russo's approximate zero-one law,
{\it Annals of Probability}, 22 (1994), 1576--1587.

\bibitem {ODo2014}
Ryan O'Donnell, Anaysis of Boolean functions, University Press, 2014.

\bibitem {Win02}  E. Winter, The Shapley value. In R.J. Aumann and S. Hart, editors, {\it Handbook of Game Theory with Economic Applications}, Vol. 3, pages 2025–2054. Elsevier Science, 2002.
\end{thebibliography}
\end{document}